\scrollmode
\documentclass[12pt]{amsart}
\usepackage{amssymb}
\usepackage{verbatim}
\usepackage{eucal}
\usepackage{eepic}
\usepackage{graphicx}
\usepackage{psfrag}
\usepackage{mathrsfs}
\swapnumbers

\theoremstyle{plain}
\newtheorem{theorem}[subsection]{Theorem}
\newtheorem{proposition}[subsection]{Proposition}
\newtheorem{lemma}[subsection]{Lemma}
\newtheorem{corollary}[subsection]{Corollary}

\theoremstyle{definition}
\newtheorem{definition}[subsection]{Definition}
\newtheorem{example}[subsection]{Example}

\theoremstyle{remark}
\newtheorem{remark}[subsection]{Remark}

\newcommand{\R}{{\mathbb R}}

\newcommand{\C}{{\mathbb C}}
\newcommand{\Z}{{\mathbb Z}}
\newcommand{\Q}{{\mathbb Q}}

\newcommand{\F}{{\mathbb F}}
\newcommand{\B}{{\mathscr B}}
\newcommand{\Proj}{{\mathbb P}}

\newcommand{\Gal}{\operatorname{Gal}}
\newcommand{\Dim}{\operatorname{dim}}
\DeclareMathOperator*{\Max}{Max}
\newcommand{\sets}[1]{[\![#1]\!]}
\newcommand{\OK}{{\mathscr{O}}}
\renewcommand{\emptyset}{\varnothing}
\renewcommand{\setminus}{\smallsetminus}

\begin{document}

\title{Symplectic Modular Symbols}

\author{Paul E. Gunnells}
\address{Department of Mathematics\\
Columbia University\\
New York, New York  10027}
\date{May 7, 1998.  Revised July 28, 1999}
\email{gunnells@math.columbia.edu}
\thanks{Partially supported by the NSF}
\begin{abstract}
Let $K/\Q $ be a number field with euclidean ring of integers $\OK$.
Let $\Gamma $ be a finite-index torsion-free subgroup of the
symplectic group $Sp_{2n}(\OK)$, and let $N$ be the cohomological
dimension of $\Gamma $.  We exhibit a finite, geometrically-defined
spanning set of $H^{N}(\Gamma ;\Z )$ by generalizing the modular
symbol algorithm of Ash and Rudolph for
$SL_{n}(\OK )$.
\end{abstract}
\maketitle

\section{Introduction}\label{introduction}
\subsection{} Let $G$ be a semisimple algebraic group defined over $\Q
$ of $\Q $-rank $\ell $, and let $X$ be the associated symmetric
space.  Let $\Gamma \subset G(\Q )$ be a torsion-free arithmetic
subgroup.  Then $H^{*}(\Gamma ;\Z )= H^{*}(\Gamma \backslash X;\Z )$,
and this cohomology vanishes for $*>N$, where $N=\Dim(X) - \ell $, the
\emph{cohomological dimension} of $\Gamma $.

The theory of modular symbols as formulated by Ash \cite{ash.minmod}
constructs an explicit spanning set for $H^{N}(\Gamma ;\Z )$ as
follows.  Let $\B$ be the Tits building associated to $G$ \cite{tits}.  By the
Solomon-Tits theorem, $\B $ has the homotopy type of a wedge of $(\ell
-1)$-spheres, and thus $\tilde{H}_{*}(\B ;\Z )$ is nonzero only in
dimension $\ell -1$.  Using the Borel-Serre compactification of the locally
symmetric space $\Gamma \backslash X$, we may
construct a map 
\begin{equation}\label{phi}
\Phi \colon H_{\ell -1}(\B ;\Z )\longrightarrow H^{N}(\Gamma ;\Z ) 
\end{equation}
that is surjective (cf. \S\ref{modular.symbols}).  Since the left-hand
side of (\ref{phi}) is generated by fundamental classes of apartments
of $\B $, this provides a geometric spanning set for $H^{N}(\Gamma )$.
These cohomology classes (or rather, their duals in homology) are
called {\em modular symbols}.

\subsection{}
The modular symbols provide a spanning set for $H^{N}(\Gamma ;\Z )$,
but they do not provide a finite spanning set, a distinction that is
important for applications.  However, suppose $K/\Q $ is a number
field with euclidean ring of integers $\OK$, and let $G(\Q ) =
SL_{n}(K)$ and $\Gamma\subset SL_{n}(\OK)$.  Then in
\cite{ash.rudolph}, Ash and Rudolph determine an explicit finite
spanning set---the {\em unimodular symbols}---and present an algorithm
to write a modular symbol as a sum of unimodular symbols
(cf. \S\ref{speciallinearcase}).  This algorithm, in conjunction with
certain explicit cell complexes, can be used to compute the action of
the Hecke operators on $H^{N}(\Gamma )$.  In turn, through work of
Ash, Pinch, and Taylor \cite{apt}, Ash and McConnell \cite{exp.ind},
and van Geemen and Top \cite{geemen.top}, this has produced
corroborative evidence for certain aspects of the ``Langlands
philosophy.''  In particular, in the case of $\Gamma \subset SL_{3}
(\Z )$, many examples of representations of the absolute Galois group
$\Gal (\bar \Q /\Q )$ have been found that appear to be associated to
cohomology classes of $\Gamma $.

\subsection{}
In this paper we solve the finiteness problem for the symplectic group:
$G(\Q ) = Sp_{2n}(K)$ and $\Gamma $ of finite index in $Sp_{2n}(\OK)$,
where $\OK$ is euclidean.  We characterize a finite spanning set
of $H^{N}(\Gamma ;\Z )$ and present an algorithm (Theorem \ref{th2})
that generalizes the modular symbol algorithm of \cite{ash.rudolph}.
To do this, we prove a relation in the homology of the
symplectic building (Theorem \ref{th1}) (see Examples
\ref{homology.sp4} and \ref{sp6.example}).

\subsection{}
We conclude this introduction by discussing the arithmetic nature of
$H^{N} (\Gamma)$ in the symplectic case, and by indicating possible
applications of Theorem \ref{th2}.

First, we consider the complex cohomology.  It is known that for
$n>1$, the groups $H^{N} (\Gamma; \C )$ do not contain cuspidal
classes (that is, cohomology classes corresponding to cuspidal
automorphic forms).  Moreover, work of Schwermer
\cite{schwermer.symp1, schwermer.symp2, schwermer.generic} shows the
following:

\begin{itemize}
\item For $n>1$, there are Eisenstein cohomology classes in the top
degree.  These classes are constructed using Eisenstein
series and characters attached to the maximal split torus of the
minimal parabolic subgroup of $G$.  For $n>3$, these are all the
Eisenstein classes appearing in the top degree.
\item Furthermore, for $Sp_{4}$, there are additional classes
constructed using Eisenstein series and cuspidal classes for $SL_{2}$.
\end{itemize}

Hence for complex coefficients these classes are either easy to
understand arithmetically, or are better studied on a lower rank group
using the classical theory of modular symbols.

Nevertheless, there is one possibility that might be of interest.  In
the case of $Sp_{6}$, one cannot exclude the possibility that there is
an Eisenstein cohomology class in the top degree associated to a
cuspidal class for $Sp_{4}$ whose infinity type is not a discrete
series representation.  It might be interesting to show the existence
of this class and to study its arithmetic nature using the results of
this paper.  (I am grateful to the referee for this suggestion.)

\subsection{}
Next, let $p$ be a prime and let $\F_{p}$ be the corresponding finite
field.  Consider the mod $p$ cohomology $H^{N} (\Gamma ;\F_{p})$, or
more generally $H^{N} (\Gamma ;M)$, where $M$ is an $\F_{p}$-module
with $\Gamma$-action.  In this setting the situation is much less
clear.  For example, there may be classes that lift to torsion classes
in integral cohomology, and thus are not associated to automorphic
forms in any obvious manner.  For discussion and examples of this, see
\cite{aac} for $\Gamma \subset GL_{3} (\Z )$, and \cite{ash.modp,
ash.sinnott} for $GL_{n}$.  One would like to know if such classes
arise in the symplectic case, and if so if the corresponding Hecke
eigenclasses are attached to Galois representations.  The algorithm in
Theorem \ref{th2}, in conjunction with the cell complex described in
\cite{mcc-macph2}, provides the means to explore this question for
$H^{4} (\Gamma )$, where $\Gamma \subset Sp_{4} (\Z )$.

\subsection{}
Finally, in another article \cite{experimental}, we describe an
algorithm to compute the Hecke action on $H^{5} (\Gamma )$, where
$\Gamma $ is a subgroup of $SL_{4} (\Z )$.  This cohomology group,
whose degree is one less than the cohomological dimension, is within
the cuspidal range.  One would like to have a symplectic version of
this algorithm that could compute the Hecke action on cuspidal classes
in $H^{3}$ of subgroups of $Sp_{4} (\Z )$ (i.e. on Siegel modular
forms of weight $3$).  The algorithm described in this paper is an
essential first step towards solving this problem.

\subsection{Acknowledgments/Related work}\label{thanks}
This paper relies heavily on the results of \cite{ash.minmod} and
\cite{ash.rudolph}.  We thank Avner Ash for much advice and
encouragement.

In \cite{mcc-macph2} the authors describe a deformation retract $W$ of
the symmetric space $Sp_{4}(\R )/U(2)$ that may be used to compute
$H^{*}(\Gamma;\Z ) $ where $\Gamma $ is a finite index subgroup of
$Sp_{4}(\Z )$.  The combinatorial data used to describe the cell
decomposition of $W$, found in \cite{mcc-macph1}, inspired the results
in Example \ref{homology.sp4} and led to Theorem \ref{th1}.  We thank
Bob MacPherson and Mark McConnell for many conversations.

Finally, we thank Mark McConnell, Richard Scott, and the referee, who
made helpful comments that improved this article.

\section{Minimal modular symbols}\label{modular.symbols}
In this section we review results on minimal modular symbols.  These
are due to Ash and Rudolph \cite{ash.rudolph} for $G=SL_{n}$ and Ash
\cite{ash.minmod} for any semisimple $\Q $-group $G$.  Our
exposition closely follows these sources.  For general results about
buildings, the reader may consult Tits \cite{tits}.

\subsection{}\label{bkground}
Let $G$, $X$, and $\Gamma $ be as in the introduction, and let
$e\in X$ be the distinguished basepoint.  Let $T$ be a maximal $\Q
$-split torus of $G$ stable under the Cartan involution corresponding
to $e$, and let $A = T(\R )^{0}$.  As $\ell $ is the $\Q $-rank of $G$,
we have that $A\cong (\R ^{>0})^{\ell }$.

Let $\bar X$ be the partial compactification of $X$ constructed by
Borel and Serre \cite{borel.serre}.  Then the closure $Z$ of $Ae$ in
$\bar X$ is homeomorphic to a closed ball of dimension $\ell $, and
$\partial Z$ lies in $\partial \bar X$.  Let $[I]\in H_{\ell
-1}(\partial \bar X)$ be the fundamental class of $\partial Z$, and
if $m\in G(\Q )$, let $[m]$ be the fundamental class of $\partial (mZ)$.

Let $\B $ be the Tits building associated to $G$.  According to
\cite[\S8.4.3]{borel.serre}, there is a homotopy equivalence $h\colon \B
\rightarrow \partial \bar X$ that takes a distinguished apartment
$A_{0}$ homeomorphically onto $\partial Z$.  This map is compatible with the
natural $G (\Q )$-action on $\B $ and $\partial \bar X$.  In
particular, if $m\in G (\Q )$, then $h (mA_{0})=\partial (mZ)$.  Since
$G (\Q )$ acts transitively on apartments, and $H_{\ell -1} (\B ;\Z )$
is generated by the fundamental classes of the apartments, we have
shown the following:

\begin{lemma}\label{gen.lemma}
\cite{ash.minmod}
The classes $[m]$ generate $H_{\ell -1}(\partial \bar X;\Z )$.

\end{lemma}

\subsection{}
Now let $\pi \colon \bar X \rightarrow \Gamma \backslash \bar X$
be the projection.  Using $\pi $ and the long exact sequence of the
pair $(\bar X, \partial \bar  X)$, we obtain
\begin{equation}\label{ex.seq}
H_{\ell -1}(\partial \bar X)\xrightarrow{\sim} H_{\ell}(\bar X,\partial \bar
X)\xrightarrow{\pi _{*}} H_{\ell}(\Gamma \backslash \bar X,\partial (\Gamma
\backslash \bar X))\xrightarrow{D} H^{N}(\Gamma \backslash \bar X) \xrightarrow{\sim}H^{N}(\Gamma \backslash X).
\end{equation}
Here the first isomorphism follows from the contractibility of $\bar
X$, the last isomorphism is induced by the canonical homotopy
equivalence $X\hookrightarrow \bar X$, and the map $D$ is Lefschetz
duality.  Since $\Gamma $ is assumed torsion-free, $D$ is an
isomorphism with integral coefficients, and $H^{N}(\Gamma \backslash
X;\Z )$ can be identified with $H^{N}(\Gamma;\Z )$.  Let
$[m]^{*}_{\Gamma }\in H^{N}(\Gamma ;\Z )$ be the image of $[m]$ under
the sequence (\ref{ex.seq}).  One of the main results of
\cite{ash.minmod} is that $\pi _{*}$ is surjective, which with Lemma
\ref{gen.lemma} implies

\begin{theorem}\label{minmod}
\cite{ash.minmod} As $m$ varies over $G(\Q )$, the classes
$[m]^{*}_{\Gamma }$ generate $H^{N}(\Gamma ;\Z )$.
\end{theorem}

\subsection{}\label{sln.bldg}
Let $K/\Q $ be a number field with ring of integers $\OK$, and let $G$ be
the algebraic $\Q $-group such that $G(\Q ) = SL_{n}(K)$.  Here $\ell
= n-1$.  We want to describe the classes $[m]$ using the combinatorics
of $\B $.  

Let $V=K^{n}$.  We assume that the basis $\{e_{1},\dots
,e_{n}\}\subset V$ are eigenvectors for the torus $T$ from
\S\ref{bkground}.  Then $\B $ is a simplicial complex with a vertex
for every proper nonzero subspace $F$ of $V$.  A set
$\{F_{1},\ldots,F_{k} \}$ of vertices of $\B $ spans a $(k-1)$-simplex
in $\B $ if and only if the corresponding subspaces can be arranged in
a proper flag:
$$0\subset F_{1} \subset \cdots \subset F_{k} \subset V. $$ 
The action of $SL_{n}(K)$ on $V$ induces an action on $\B $ and on
$H_{*}(\B ;\Z )$.  

To construct the classes $[m]$, we use an auxiliary simplicial
complex.  Let $\sets{n}$ be the set $\{1,\ldots,n \}$.  Let $\partial
\Delta _{n-1}$ be the barycentric subdivision of the boundary complex
of the $(n-1)$-simplex.  In other words, $\partial \Delta _{n-1}$ is a
simplicial complex with vertices corresponding to the proper non-empty
subsets $I$ of $\sets{n}$, and a collection of subsets
$I_{1},\ldots,I_{k}$ corresponds to a $(k-1)$-simplex in $\partial
\Delta _{n-1}$ if and only if they can be arranged in a proper flag:
$I_{1}\subset \cdots \subset I_{k}$.  We may orient $\partial \Delta
_{n-1}$ using the standard ordering on $\sets{n}$.

Following \cite{ash.rudolph}, given $n$ points in $V\setminus\{0 \}$, we may
construct a class in $H_{n-2}(\B ;\Z )$ as follows.  Given
$v_{1},\ldots,v_{n}\in V\setminus\{0 \}$, we define a simplicial map
\begin{equation}\label{delta}
\phi \colon \partial \Delta _{n-1}\longrightarrow \B 
\end{equation}
by sending the vertex $I\subset \sets{n}$ to the vertex of $\B $
corresponding to the subspace spanned by $\left\{ v_{i} \,|\, i\in I
\right\}$.  If $\xi $ is the fundamental class of the geometric
realization of $\partial \Delta _{n-1}$, then $\phi _{*}(\xi )$ is a
class in $H_{n-2}(\B ;\Z )$.  Thus under the composition 
\[
\partial \Delta _{n-1}\xrightarrow{\phi }\B \xrightarrow{h}\partial
\bar X,
\]
we have constructed a class in $H_{n-2} (\partial \bar X; \Z )$.

\begin{definition}\label{ms}
\cite{ash.rudolph} The {\em modular symbol} associated to
$v_{1},\ldots,v_{n}\in V\setminus\{0 \}$ is the class in $H_{n-2}
(\partial \bar X; \Z )$ constructed above, and is denoted
$[v_{1},\ldots,v_{n}]$.  If $m\in M_{n}(K)$ ($n\times n$ matrices over
$K$), then by $[m]$ we mean the modular symbol constructed using the
columns of $m$.
\end{definition}

\begin{proposition}\label{equivalence.sl}
If $m\in SL_{n} (K)$, then the construction of $[m]$ given in
Definition \ref{ms} coincides with that given in \S\ref{bkground}.  In
particular, the classes $[m]$ span $H_{n-2} (\partial \bar X; \Z )$.
\end{proposition}

\begin{proof}
Recall that we have a homotopy equivalence $h\colon \B \rightarrow
\partial \bar X $ taking a distinguished apartment $A_{0}$
homeomorphically onto $\partial Z$ (\S\ref{bkground}).  This apartment
is in fact $\phi (\partial
\Delta _{n-2}) $, where $\partial
\Delta _{n-2}$ is constructed using the basis $e_{1},\dots ,e_{n}\in
V$.  These identifications are compatible with the action of $SL_{n}
(K)$ on $V$, $\B $, and $\partial \bar X $, so the result follows.
\end{proof}

Modular symbols have the following properties.

\begin{proposition}\label{props}
\cite{ash.rudolph} The map $[\phantom{m}]\colon
M_{n}(K)\rightarrow H_{n-2}(\partial \bar X;\Z )$ satisfies the following:
\begin{enumerate}
\item $[v_{1},\ldots,v_{n}] = (-1)^{|\tau |}[\tau (v_{1}),\ldots,\tau (v_{n})]$, where $\tau \in S_{n}$ is a
permutation on $n$ letters, and $|\tau |$ is the length of $\tau $.
\item If $q\in K$, then $[q v_{1},v_{2},\ldots,v_{n}] = [v_{1},\ldots,v_{n}]$.
\item If the $v_{i}$ are linearly dependent, then
$[v_{1},\ldots,v_{n}] = 0$.
\item If $v_{0},\ldots,v_{n}\in V\setminus\{0 \}$, then 
$$\sum _{i}(-1)^{i}[v_{0},\ldots,\hat{v_{i}},\ldots,v_{n}]=0. $$
\end{enumerate}
Furthermore, $[\phantom{m}]$ is surjective.
        
\end{proposition}

\subsection{}\label{speciallinearcase}
Now assume that $\OK$ is a euclidean ring with respect to a
multiplicative norm $\|\phantom{m}\|\colon \OK\rightarrow \Z ^{\geq
0}$.  Using multiplicitivity extend the norm to a map
$\|\phantom{m}\|\colon K\rightarrow \Q ^{\geq 0}$.  We
recall how to identify a $\Gamma $-finite spanning set of modular
symbols.  By a {\em primitive} vector, we mean a vector $v\in \OK^{n}$
such that the greatest common divisor of the entries of $v$ is a unit.

Let $L$ be an $\OK$-submodule of $\OK^{n}$ of rank $k\leq n$.  Since
$\OK $ is a principal ideal domain, $L$ has a free $\OK $-basis $B =
\{v_{1},\dots ,v_{k}\}$.  Choose $W=\{w_{k+1},\dots ,w_{n} \}\subset
\OK^{n} $ such that $B\bigsqcup W$ is a $K$-basis of $K^{n}$, and $W$
may be extended to an $\OK $-basis of $\OK ^{n}$.

We define the {\em index} of $L$ by 
\[
i(L) := \| \det (v_{1},\ldots,v_{k},w_{k+1},\dots w_{n}) \|.
\]
It is easy to see that $i (L)$ is independent of the choices above,
and that 
\[
i (L)=1 \quad \hbox{if and only if} \quad (L\otimes _{\OK }K)\cap \OK ^{n} = L.
\]
Furthermore, if $L$ has rank $n$ and $\|\phantom{m}\|$ is the usual
norm $N_{K/\Q }\colon K\rightarrow \Q $, then $i (L) = [\OK ^{n}: L]$.

We write $i(v_{1},\ldots,v_{k})$ for $i(L)$ if we are given a specific
set of linearly independent vectors generating $L$.  If the $v_{i}$
are linearly dependent, we define $i(v_{1},\ldots,v_{k})=0$.

\begin{definition}\label{cand}
Let $v_{1},\ldots,v_{k}\in \OK^{n}$ be linearly independent primitive
vectors.  Then a {\em candidate} for the $v_{i}$ is a
primitive $x\in \OK^{n}\smallsetminus \{0 \}$ satisfying
\[
0\leq i(x,v_{1},\ldots,\hat v_{i},\ldots,v_{k}) <
i(v_{1},\ldots,v_{k}), \quad \hbox{$1\leq i\leq k$.} 
\]
\end{definition}

A fundamental result is the following:

\begin{proposition}\label{cand.exist}
\cite{ash.rudolph} Let $v_{1},\ldots,v_{k}$ be a linearly independent
set of primitive vectors.  If $i(v_{1},\ldots,v_{k})>1$, then a
candidate $x$ for the $v_{i}$ exists.
\end{proposition}

\begin{proof}
We give the proof since our statement differs slightly from that
found in \cite{ash.rudolph}.  First assume $k=n$.  Let $L$ be the
lattice spanned by the $v_{i}$, and let $w$ be a primitive vector in
$\OK ^{n}$ that is not in $L$.  Note that $w$ exists since $i(L)>1$.

Let $A\in M_{n}(\OK )$ be the matrix with columns
$v_{1},\ldots,v_{n}$, and let $A_{i}[w]$ be the matrix obtained by
replacing the $i$th column of $A$ with $w$.  Since $\OK $ is
euclidean, there exist $\alpha _{i},\beta _{i}\in \OK $ such that
$$\det A_{i}[w]=\alpha _{i}\det A + \beta_{i}\quad \hbox{where $0\leq \| \beta _{i}
\| < \| \det A \|$}. $$

Now let $x=w-\sum_{i} \alpha _{i}v_{i}$.  By our choice of $w$ the
vector $x$ is nonzero.  It is easy to check that
$\det A_{i}[x] = \beta _{i}$.  Since $0\leq \| \beta _{i}
\| < \| \det A \|$ and $\|\det A_{i}[x]\| = i(x,v_{1},\ldots,\hat
v_{i},\ldots,v_{n})$, the result follows.

Now assume $k<n$.  Since $\OK $ is euclidean, $L$ is a free module.
Hence we can choose an isomorphism $L\otimes K\rightarrow K^{k}$ that
carries $(L\otimes K)\cap \OK ^{n}$ onto $\OK ^{k}$. This
brings us back to the full rank setting, and we may argue as
above.
\end{proof}

Notice that $x$ can be written as $\sum q_{i}v_{i}$, with $q_{i}\in K$
satisfying $0\leq \|q_{i}\| < 1$.  Furthermore, since the $v_{i}$ are
linearly independent, at least one $q_{i}\not =0$.

\begin{theorem}\label{finite}
\cite{ash.rudolph} As $m$ ranges over $SL_{n}(\OK)$, the classes $[m]$
generate $H_{n-2}(\partial \bar X;\Z )$.  Hence if $\Gamma \subset
SL_{n}(\OK)$ is torsion-free and of finite index, then the classes
$[m]^{*}_\Gamma $ provide a finite spanning set of $H^{N}(\Gamma ;\Z
)$.
\end{theorem}

\begin{proof}
Repeatedly applying Proposition \ref{cand.exist} and Proposition
\ref{props} (4), we may write any class $[m]$ as a sum $\sum [m_{i}]$,
where the determinant of each $m_{i}$ is a unit.  Then applying
Proposition \ref{props} (2), we may divide the first column of $m_{i}$
by the determinant of $\det m_{i}$ to get $m'_{i}\in SL_{n} (\OK )$
satisfying $[m_{i}]=[m'_{i}]$.
\end{proof}

\begin{definition}\label{unimod}
The classes $\bigl \{[m] \bigm |m\in SL_{n}(\OK) \bigr\}$ are called
{\em unimodular symbols}.
\end{definition}

\section{Symplectic modular symbols}\label{symp.building}
In this section we generalize the results of \S\ref{sln.bldg} to the
symplectic case.  Sections \ref{notation} and \ref{sp2n.bldg} recall
well-known facts about symplectic geometry and the building associated
to the symplectic group.  In \S\ref{symp.mod.syms} we translate
results of \S\ref{sln.bldg} to the symplectic setting, and in
\S\S\ref{part4}-\ref{the.proof} we describe a relation in the homology
of the building that is crucial to our finiteness result.  Throughout
this section we do not assume that the ground field $K$ has a
euclidean ring of integers.

\subsection{}\label{notation}
First we recall some elementary facts about symplectic geometry to fix
notation.  Fix a field $K$, and let $V$ be the vector space $K^{2n}$.
If $k\in \Z $, we use the notation $\bar k$ for $2n+1-k$.  We fix a
non-degenerate alternating bilinear form
$\langle\phantom{l},\phantom{l}\rangle\colon V\rightarrow K$.  A
basis $v_{1},\ldots,v_{n},v_{\bar n},\ldots,v_{\bar 1}$ of $V$ is said
to be {\em symplectic} if
\[
\langle v_{i}, v_{j}\rangle=\begin{cases}
        1 & \text{if $j=\bar \imath$ and $i<j$},\\
        -1 & \text{if $j=\bar \imath$ and $i>j$},\\
        0 & \text{otherwise}.
\end{cases}
\]
We let $G=Sp_{2n}(K)$ be the subgroup of $GL_{2n}(K)$ preserving the form.
Thus $Sp_{2n}(K) = \left\{g\in GL_{2n}(K) \mid \langle gv, gw\rangle =
\langle v, w \rangle \right\}$.  The group $G$ has $\Q $-rank $\ell =n$.

Given any $x\in V$, we define $x^{\perp}$ to be the set $\left\{y\in
V\mid\langle x, y\rangle = 0 \right\}$.  If $x$ is nonzero, then $x^{\perp}$ is
a hyperplane containing $x$.  A subspace $F\subset V$ is called {\em
isotropic} if $\langle v,w \rangle =0$ for all $v,w\in F$.  Any
one-dimensional subspace is isotropic, and the largest dimension an
isotropic subspace may have is $n$.  An $n$-dimensional isotropic
subspace is called a {\em Lagrangian subspace}.

\subsection{}\label{sp2n.bldg}
From now on we exclusively use $\B $ to denote the building associated
to $Sp_{2n}(K)$.  We have the following description of $\B $ as a
simplicial complex, analogous to that found in \S\ref{sln.bldg}:
vertices of $\B $ correspond to nonzero isotropic subspaces of $V$,
and simplices of $\B $ correspond to flags of nonzero isotropic
subspaces.  

To describe the geometry of $\B $ we must use
cross-polytopes instead of simplices, and so we recall their
definition.  Let $e_{1},\ldots,e_{n}$ be the standard basis of $\R
^{n}$.  Then a {\em cross-polytope} on $2n$ vertices is a polytope
isomorphic to the convex hull of the points $\pm
e_{1},\ldots,\pm e_{n}$.  Examples are the square ($n=2$) and the
octahedron ($n=3$).

Since the proper faces of cross-polytopes are simplicial, we may
regard their boundary complexes as simplicial complexes.  Let
$\partial \beta _{n}$ be the first barycentric subdivision of the
boundary complex of the cross-polytope on $2n$ vertices.  To describe
the structure of $\partial \beta _{n}$, we use the notation
$\sets{n}^{\pm } := \{1,\ldots,n,\bar n,\ldots,\bar 1 \}$.  We order
this set by $1<\dots <n<\bar n <\dots <\bar 1$.

\begin{definition}\label{isotropic.subsets}
(cf. \cite{gel.serg}) A nonempty subset $I\subset \sets{n}^{\pm }$ is
called {\em isotropic} if for all $i,j\in I$, we have $i\not =\bar
\jmath$.
\end{definition} 
Note that $\# I \leq n$ for any isotropic $I$.  As in the $SL_{n}$
case, vertices of $\partial \beta _{n}$ correspond to isotropic
subsets of $\sets{n}^{\pm }$: if we take $\beta $ to be the convex
hull of the points $\pm e_{1},\ldots,\pm e_{n}$, then the set $I$
corresponds to the linear span of $\{e_{i}\mid i\in I \}$, where
$e_{\bar \imath } := -e_{i}$.  Similarly, simplices of $\partial \beta
_{n}$ correspond to proper flags of isotropic subsets.

\subsection{}\label{symp.mod.syms}
Now suppose that we are given $2n$ points $v_{1},\ldots,v_{n},v_{\bar
n},\ldots,v_{\bar 1}\in V\setminus\{0 \}$.  
\begin{definition}\label{isoptropycondition}
The set $v_{1},\ldots,v_{n},v_{\bar n},\ldots,v_{\bar 1}$ is said to
satisfy the {\em isotropy condition} if for every isotropic subset
$I\subset \sets{n}^{\pm }$, the subspace spanned by $\{v_{i}\,|\,i\in
I \}$ is isotropic.
\end{definition}

Notice that this condition is strictly weaker than requiring that the
$v_{i}$ form a symplectic basis.  In particular, the columns of any
$m\in Sp_{2n}(K)$ satisfy this condition.  Following \S\ref{sln.bldg},
we define a simplicial map
\[
\phi \colon \partial \beta _{n}\longrightarrow \B 
\]
by sending the vertex corresponding to an isotropic $I\subset
\sets{n}^{\pm }$ to the vertex of $\B $ corresponding to the linear
span of $\{v_{i}\,|\,i\in I \}$.  Via the composition $\phi \circ h$,
we have an induced map on homology taking the fundamental class $\xi $
of the geometric realization of $\partial \beta _{n}$ to a class in
$H_{n-1}(\partial \bar X;\Z )$.

\begin{definition}\label{sms.def}
Given $v_{1},\ldots,v_{n},v_{\bar n},\ldots,v_{\bar 1}\in
V\setminus\{0 \}$ satisfying the isotropy condition, let
$[v_{1},\ldots,v_{n};v_{\bar n},\ldots,v_{\bar 1}]\in H_{n-1}(\partial
\bar X;\Z )$ denote the class constructed above.  This class is called
a {\em symplectic modular symbol}.  If the columns of $m\in M_{2n}(K)$
satisfy the isotropy condition, we denote the symplectic modular
symbol corresponding to its columns by $[m]$.
\end{definition}

Symplectic modular symbols satisfy properties similar to those
satisfied by the special linear symbols.  For instance, minor
modification of the proof of Proposition \ref{equivalence.sl} shows
that the construction of $[m]$ in Definition \ref{sms.def} agrees with
that of \S\ref{bkground}, and thus the modular symbols span
$H_{n-1}(\partial \bar X;\Z )$.  Furthermore, we have the following
analog of the first parts of Proposition \ref{props}, whose proof is
easily checked:

\begin{proposition}\label{sms.props}
Symplectic modular symbols enjoy the following properties:
\begin{enumerate}
\item [1a.] Let $\tau \in S_{n}$ be a permutation on $n$ letters.  Then
$$[v_{1},\ldots, v_{n};v_{\bar n},\ldots,v_{\bar 1}] = [\tau
(v_{1}),\ldots, \tau (v_{n});\tau (v_{\bar n}),\ldots,\tau (v_{\bar
1})],
$$  
where $\tau(v_{k}) := v_{\tau (k)}$ and $\tau (v_{\bar k}) :=
v_{\overline{\tau (k)}}$ for $k\in \sets{n}$.
\item [1b.] $[v_{1},v_{2},\ldots,v_{n};v_{\bar n},\ldots,v_{\bar 2},v_{\bar 1}] = 
-[v_{\bar 1},v_{2},\ldots,v_{n};v_{\bar n},\ldots,v_{\bar 2},v_{1}]$.
\item [2.] If $q\in K$, then $[q v_{1},v_{2},\ldots,v_{n};v_{\bar n},\ldots,v_{\bar 1}] = [v_{1},\ldots,v_{n};v_{\bar n},\ldots,v_{\bar 1}]$.
\item [3.] If the $v_{i}$ are linearly dependent, then
$[v_{1},\ldots,v_{n};v_{\bar n},\ldots,v_{\bar 1}] = 0$.
\end{enumerate}
\end{proposition}

\subsection{}\label{examples}
The symplectic analogue of Proposition \ref{props} (4) is more
complicated and forms one of the main results of
this article.  To motivate the result we illustrate it for $Sp_{4}(K)$
and $Sp_{6}(K)$.  In light of the homotopy equivalence $\B \rightarrow
\partial \bar X$ (\S\ref{bkground}), we work in $H_{*} (\B )$
rather than $H_{*} (\partial \bar X)$.

\begin{example}\label{homology.sp4}
Let $G = Sp_{4}(K)$.  Then the simplicial complex
$\B $ is a graph with two types of vertices, corresponding to the one-
and two-dimensional isotropic subspaces of $V$.

Let $[m]=[v_{1},v_{2};v_{\bar 2},v_{\bar 1}]$ be a symplectic modular
symbol for $Sp_{4}(K)$.  We may picture the subspace configuration in
$V=K^{4}$ determined by $\{v_{i}\}$ by passing to $\Proj ^{3}(K)$.
Figure \ref{config} shows the projectivized configuration on the left,
along with the apartment in $\B $ corresponding to $[m]$ on the right.
In $\B $ we represent the one-dimensional (respectively
two-dimensional) isotropic subspaces of $V$ by solid (resp. hollow)
vertices.  By abuse of notation we use the same symbol for a point in
$K^{4}$, the point it determines in $\Proj ^{3}(K)$, and the vertex it
determines in $\B $.  The lines on the left of Figure \ref{config} are
the projectivizations of the Lagrangian planes determined by the
$v_{i}$; we have not drawn the images of the two non-Lagrangian
planes.

\begin{figure}[h]
\psfrag{v1}{$v_{1}$}
\psfrag{v2}{$v_{2}$}
\psfrag{vb1}{$v_{\bar 1}$}
\psfrag{vb2}{$v_{\bar 2}$}
\centerline{\includegraphics[scale = .7]{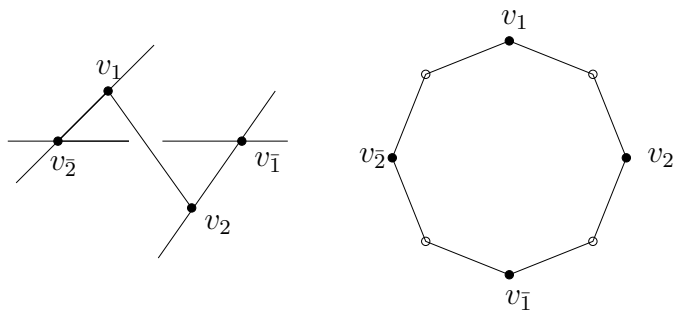}}
\caption{\label{config} A configuration in $\Proj ^{3}(K)$ and its corresponding apartment.}
\end{figure}

Now choose $x\in V\smallsetminus \{0 \}$.  We will use $x$ to construct a
class in $H_{1}(\B ;\Z )$ homologous to $[m]$.  Recall that
$x^{\perp}$ is the set of all $y\in V$ satisfying $\langle x, y
\rangle = 0$.  In $\Proj ^{3}(K)$, $x^{\perp}$ becomes a hyperplane
that, for generic $x$, determines four new points by intersection
with the original Lagrangian lines (see Figure \ref{newpoints}).
\begin{figure}[h]
\begin{center}
\psfrag{v1}{$v_{1}$}
\psfrag{v2}{$v_{2}$}
\psfrag{vb1}{$v_{\bar 1}$}
\psfrag{vb2}{$v_{\bar 2}$}
\psfrag{x}{$x$}
\includegraphics[scale = .7]{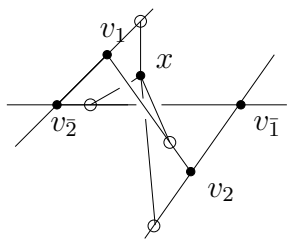}
\end{center}
\caption{\label{newpoints}Constructing new points.}
\end{figure}
These new points and lines determine four new apartments in $\B $, and
provide a relation $[m] = \sum_{i=1}^{4} [m_{i}]$ in homology, as in Figure
\ref{hom.fig}. 
\begin{figure}[h]
\begin{center}
\psfrag{v1}{$v_{1}$}
\psfrag{v2}{$v_{2}$}
\psfrag{vb1}{$v_{\bar 1}$}
\psfrag{vb2}{$v_{\bar 2}$}
\psfrag{x}{$x$}
\psfrag{a}{$a$}
\psfrag{b}{$b$}
\psfrag{c}{$c$}
\psfrag{d}{$d$}
\includegraphics[scale = .7]{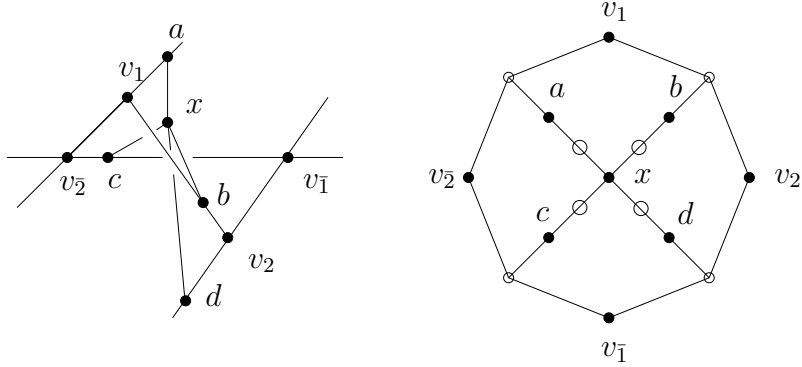}
\end{center}
\caption{\label{hom.fig}A homology relation.} 
\end{figure}

Now suppose that $x$ is not generic.  Then for some subset $I\subset
\sets{2}^{\pm }$, we have $\langle x, v_{i}\rangle = 0$ for $i\in I$.
Up to symmetry there are three nontrivial possibilities, which we depict in
Figure \ref{nongensq.fig}.
The corresponding apartments appear in Figure \ref{nongenapt.fig}.
Note that in these cases there are fewer apartments in the images of
the relations.
\begin{figure}[h]
\begin{center}
\psfrag{v1}{$v_{1}$}
\psfrag{v2}{$v_{2}$}
\psfrag{vb1}{$v_{\bar 1}$}
\psfrag{vb2}{$v_{\bar 2}$}
\psfrag{x}{$x$}
\includegraphics[scale = .7]{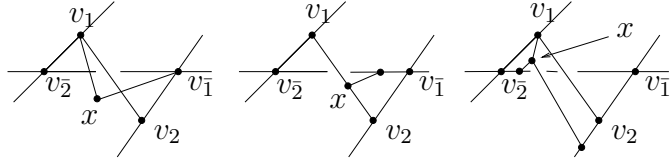}
\end{center}
\caption{\label{nongensq.fig}$I = \{1,\bar 1 \}$,$\{1,2 \}$, and $\{1 \}$.}
\end{figure}
\begin{figure}[h]
\begin{center}
\psfrag{v1}{$v_{1}$}
\psfrag{v2}{$v_{2}$}
\psfrag{vb1}{$v_{\bar 1}$}
\psfrag{vb2}{$v_{\bar 2}$}
\psfrag{x}{$x$}
\includegraphics[scale = .7]{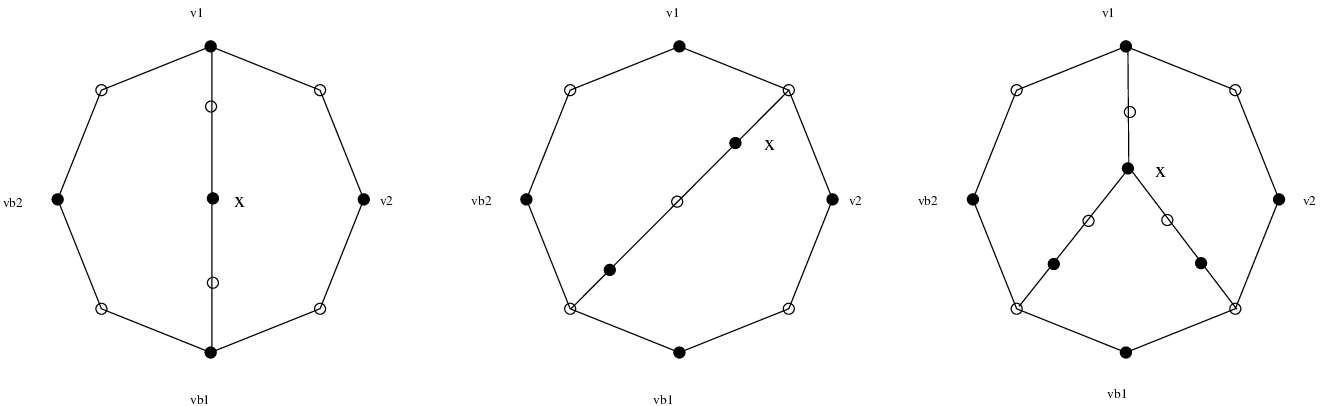}
\end{center}
\caption{\label{nongenapt.fig}}
\end{figure}
\end{example}

\begin{example}\label{sp6.example}
Consider the case $G=Sp_{6}(K)$.  Now $\B $ has three kinds of
vertices, corresponding to isotropic points, lines, and planes in
$\Proj^{5}(K)$.  The configuration in $\Proj^{5}(K)$ corresponding to
a symplectic modular symbol $m$ is combinatorially an octahedron with
isotropic faces.  Choosing a generic point $x$, we construct twelve
new points by intersecting $x^{\perp}$ with the isotropic lines of the
octahedron.  In Figure \ref{octa.fig} we show the configuration
corresponding to $[m]$ along with these constructed points, which are
shown as hollow dots.  (Although the configuration properly lives in
$\Proj^{5}(K)$, we show it in three dimensions for clarity.) Notice
that the constructed points satisfy nontrivial linear dependencies:
three constructed points are collinear if they lie on an isotropic
plane corresponding to a facet of the original octahedron.  (These
dependencies are only true in $\Proj^{5}(K)$.)  These dependencies
ensure the relation $[m] = \sum [m_{i}]$ in homology.

\begin{figure}[h]
\begin{center}
\includegraphics[scale = .7]{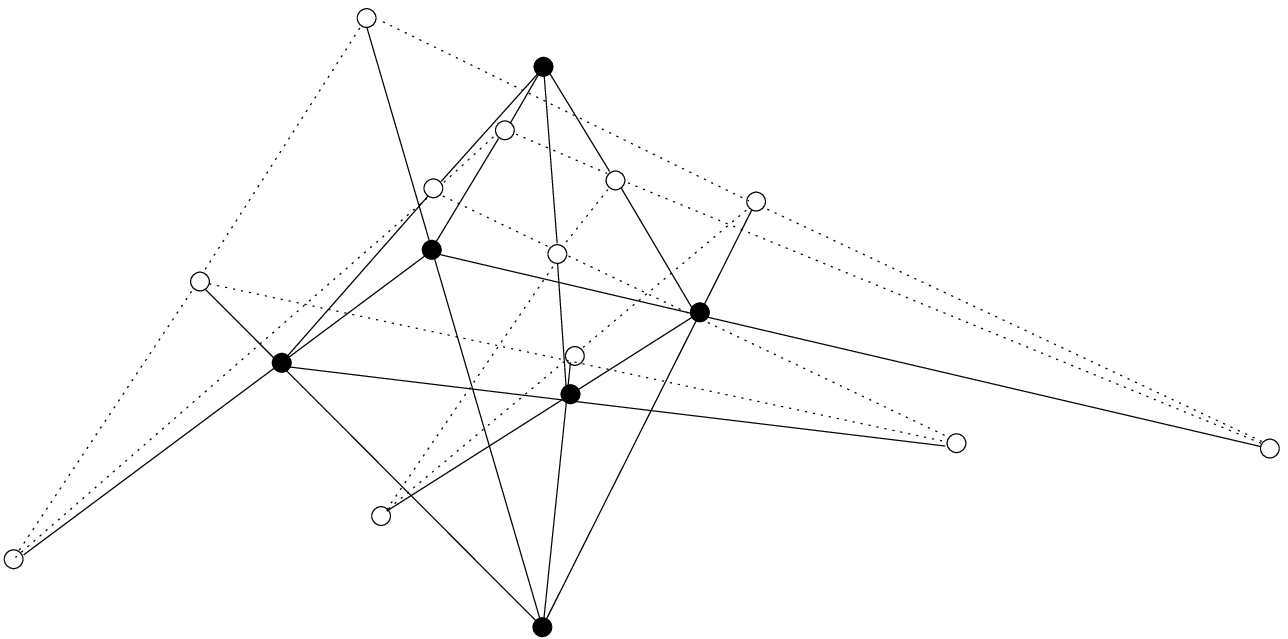}
\end{center}
\caption{\label{octa.fig}A configuration in $\Proj^{5}(K)$.}
\end{figure}
\end{example}

\subsection{}\label{part4}
Returning to the general case, let $[m]=[v_{1},\ldots,v_{n};v_{\bar
n},\ldots,v_{\bar 1}]$ be a symplectic modular symbol.  Let $x\in
V\setminus\{0 \}$, and let $x^{\perp}$ be as above.  Let $D_{x}\subset
\sets{n}^{\pm }$ be the set of indices such that $\langle
x,v_{i}\rangle = 0$ if and only if $i\in D_{x}$.  Given distinct $i,j\in
\sets{n}^\pm $ with $i\not =\bar \jmath$ and not both $i,j\in D_x$,
define $x_{ij}$ by
\begin{equation}\label{xij}
x_{ij} = \langle x,v_{i}\rangle v_{j} - \langle x,v_{j}\rangle v_{i}.
\end{equation}
If both $i$ and $j$ lie in $D_x$, then $x_{ij}$ is not defined.  The
point $x_{ij}$ lies on the intersection of $x^{\perp}$ with the
isotropic plane spanned by $v_{i}$ and $v_{j}$.

Now we define new matrices built from $m$, $x$, and the $x_{ij}$.

\begin{definition}\label{msubi}
Let $[m]$ be a symplectic modular symbol.  Choose $x\in
V\smallsetminus \{0 \}$, and construct the $x_{ij}$ as in \eqref{xij}.
If $i\not \in D_{x}$, define the matrix $m_{i}$ to be the matrix
obtained by altering $m$ according to the following rules:
\begin{enumerate}
\item Replace $v_{\bar \imath}$ by $x$.
\item For $j\in \sets{n}^{\pm }\setminus\{i,\bar \imath \}$, 
replace $v_{j}$ by $x_{ij}$.
\end{enumerate}
\end{definition}
A priori $[m_{i}]$ may not be a symplectic modular symbol, since
its columns might not satisfy the isotropy condition.  Hence we
state 
\begin{proposition}
For each $i\in \sets{n}^{\pm }\setminus D_{x}$, each $[m_{i}]$ is a symplectic
modular symbol.
\end{proposition}

\begin{proof}
It is easy to check that $\langle x_{ij},x_{ik}\rangle = 0$ and
$\langle v_{i}, x_{ij} \rangle =0$ in all the necessary cases.  Since
$\langle x, x_{ij} \rangle =0$ by construction, the result follows.
\end{proof}

As indicated in Example \ref{sp6.example}, the $x_{ij}$ satisfy linear
dependencies, which we record in the following lemma.

\begin{lemma}\label{linearspan}
Suppose $I=\{i,j,k \}$ is an isotropic subset, and $\# ( I\cap D_{x})\leq
1$.  Then the points $x_{ij}$, $x_{jk}$, and $x_{ik}$ are linearly
dependent.  
\end{lemma}

\begin{proof}
The points satisfy the identity
\[
\langle x, v_{k} \rangle x_{ij }= \langle x, v_{j} \rangle
x_{ik}-\langle x, v_{i} \rangle x_{jk}.
\]
\end{proof}

\subsection{}\label{the.proof}
We come now to the main result of this section.

\begin{theorem}\label{th1}
Let $[m]=[v_{1},\ldots,v_{n};v_{\bar n},\ldots,v_{\bar 1}]$ be a
symplectic modular symbol for $Sp_{2n} (K)$, and choose $x\in
V\smallsetminus \{0 \}$.  Define $D_{x}$ as above, and let $[m_{i}]$
be the symplectic modular symbols constructed in Definition
\ref{msubi}.  Then in $H_{n-1}(\partial \bar X;\Z )$, we have
\begin{equation}\label{homology}
[m]=\sum_{i\in \sets{n}^{\pm }\setminus D_{x}} [m_{i}]. 
\end{equation}
\end{theorem}

\begin{proof}
As in Examples \ref{homology.sp4} and \ref{sp6.example}, we prove the
relation in $H_{n-1} (\B ;\Z )$.
Let $A$ (respectively $A_{i}$) be the apartment corresponding to the
symplectic modular symbol $[m]$ (respectively $[m_{i}]$).  We
think of these apartments as being explicit simplicial cycles in $\B
$, and we will show that $[m]=\sum
_{i}[m_{i}]$ by examining these cycles.  

We begin by fixing some notation.  Let $(a_{1},\ldots,a_{k})$ be
an ordered tuple of linearly independent points of $V$ lying in a
Lagrangian subspace, and let $F(a_{1},\dots ,a_{k})\subset V$ be
their linear span.  Let $\sigma (a_{1},\ldots,a_{k})\subset \B $
be the simplex corresponding to the flag
\[
0\subset F(a_{1})\subset F(a_{1},a_{2})\subset \cdots \subset F(a_{1},\dots ,a_{k})\subset V.
\]
Recall that the {\em (closed) star}
of a simplex $\sigma $ in a simplicial complex is the set of all
simplices $\sigma '$ meeting $\sigma $, as well as the faces of all
such $\sigma '$.  Also recall that maximal simplices in $\B $ are
called {\em chambers}.

First assume that the point $x$ is generic with respect to the
$v_{i}$, so that $D_{x}=\emptyset $.
Consider the column vector $v_{i}$ from $[m]$.  The chambers in
$A$ appearing in the star of $v_{i}$ are the simplices of the form 
\begin{equation}\label{sigma}
\sigma (v_{i},v_{k_{1}},\ldots,v_{k_{n-1}})\subset A,
\end{equation}
where $\{i, k_{1},\ldots,k_{n-1} \}\subset
\sets{n}^{\pm }$ is isotropic.  
These chambers correspond to the flags 
\begin{equation}\label{onef}
0\subset F(v_{i})\subset F(v_{i},v_{k_{1}})\subset \cdots \subset F(v_{i},v_{k_{1}},\ldots,v_{k_{n-1}})\subset V. 
\end{equation}
On the right of \eqref{homology}, in $A_{i}$, we have the chambers
\begin{equation}
\sigma (v_{i},x_{i,k_{1}},\ldots,x_{i,k_{n-1}})\subset A_{i}.
\end{equation}
These chambers correspond to the flags
\begin{equation}\label{twof}
0\subset F(v_{i})\subset F(v_{i},x_{i,k_{1}})\subset \cdots \subset F(v_{i},x_{i,k_{1}},\ldots,x_{i,k_{n-1}})\subset V.
\end{equation}
The sets of flags in (\ref{onef}) and (\ref{twof}) coincide by
the definition of the $x_{ij}$, and these chambers appear with the
same orientations on both sides of (\ref{homology}). Taking all
permutations of $\{k_{1},\dots k_{n-1} \}$ in
\eqref{sigma}--\eqref{twof}, we obtain all chambers in the star of
$v_{i}$.  Hence any chamber in $A$ appears once in a unique $A_{i}$
with the same orientation.

Now we claim that each of the remaining chambers in the $A_{i}$
appears exactly twice with opposite orientations.  Any such chamber
must appear in the star of $x$ or $x_{ij}$, and we first consider the
star of $x$.  Choose an apartment $A_{i}$, a point $x_{ij}$, and
let $I\subset \sets{n}^{\pm }$ be a maximal isotropic subset of the form
$\{i,j,k_{1},\dots k_{n-2}\}$.  Let $F_{I}$ be the Lagrangian subspace
corresponding to $I$.  Consider the chambers with $x$ and $x_{ij}$ as
vertices, and with all vertices other than $x$ corresponding to
subspaces lying in $F_{I}$.  These chambers have the form
\begin{equation}\label{one}
\sigma (x,x_{ij},x_{i,k_{1}},\ldots,x_{i,k_{n-2}}) \subset A_{i}
\end{equation}
or
\begin{equation}\label{two}
\sigma (x,x_{ij},x_{j,k_{1}},\ldots,x_{j,k_{n-2}}) \subset A_{j}.
\end{equation}
In \eqref{one} and \eqref{two}, we allow all permutations of
$\{k_{1},\dots k_{n-2} \}$.  By Lemma \ref{linearspan}, these chambers
correspond to the same isotropic flag.  Furthermore, it is not
difficult to see that the chambers in (\ref{one}) and (\ref{two})
appear in $A_{i}$ and $A_{j}$ with opposite orientations.  We may
apply this argument to any pair $\{i,j \}$ and any $F_{I}$ with $\{i,j
\}\subset I$, and so all the chambers in the star of $x$ in the
$A_{i}$ cancel each other in \eqref{homology}.

To complete the proof for generic $x$, we investigate any remaining
chambers on the right-hand side of \eqref{homology}.  These must
appear in the stars of the $x_{ij}$.  Fix $x_{ij}$, and let
$I=\{i,j,k_{1},\dots k_{n-2}\}$ and $F_{I}$ be as above.  The chambers
meeting the star of $x_{ij}$ and with all vertices
corresponding to subspaces lying in $F_{I}$ are of two types: first,
\begin{align}
\sigma(x_{ij},x_{i,k_{1}},\ldots,x_{i,k_{n-2}},v_{i})
&\subset A_{i}\quad \hbox{and}\label{one12:a} \\
\sigma(x_{ij},x_{j,k_{1}},\ldots,x_{j,k_{n-2}},v_{j})
&\subset A_{j}\label{one12:b},\\
\intertext{and secondly}
\sigma(x_{ij},x_{i,k_{1}},\ldots,x_{i,k_{n-2}},x) &\subset A_{i}\quad \hbox{and}\label{one12:c}\\
\sigma (x_{ij},x_{j,k_{1}},\ldots,x_{j,k_{n-2}},x) &\subset A_{j}\label{one12:d}.
\end{align} 
In \eqref{one12:a}--\eqref{one12:d}, we allow all permutations of the
right $n-1$ vertices.  By Lemma \ref{linearspan}, the isotropic flags
corresponding to \eqref{one12:a} and \eqref{one12:b} (respectively
\eqref{one12:c} and \eqref{one12:d}) coincide, and checking
orientations shows that these cancel in pairs.  This accounts for all
the chambers on both sides of (\ref{homology}), and hence the result
follows for generic $x$.

Now assume that $D_{x}\not =\emptyset $.  Write $D_{x}=I\bigsqcup J$,
where $I=\bar I$ and $J\cap \bar J=\emptyset $.  We claim it is
sufficient to assume $I=\emptyset $.  Indeed, let
$S=\sets{n}^{\pm }\setminus I$, and let $V'\subset V$ be the span of
$\{v_{i}\mid i\in S \}$.  Then $x\in V'$, and $V'$ is a symplectic
space whose form $\langle\phantom{a},\phantom{a}\rangle'$ is the
restriction of $\langle\phantom{a},\phantom{a}\rangle$.  Furthermore,
the vectors $\{v_{i}\mid i\in S \}$ define an apartment in the
building $\B '$ associated to
$(V',\langle\phantom{a},\phantom{a}\rangle')$, and thus determine a
symplectic modular symbol $[m']\in H_{n'-1} (\B ';\Z )$, where $n'=\frac{1}{2}\#S$.  Writing $[m]=\sum [m_{i}]$ in $H_{n-1} (\B;\Z )$ is equivalent
to writing $[m']=\sum [m_{i}']$ in $H_{n'-1} (\B;\Z )$.  Hence, by
induction we may assume that $D_{x}$ contains no subset $I$ with $\bar
I = I$.

This is the same as $D_{x}$ being isotropic, and so up to symmetry the
only invariant of $D_{x}$ is its cardinality.  As before we proceed by
investigating the stars of vertices.  We merely indicate which
chambers cancel which and omit the details.

We first consider the case that $\# D_{x}\leq n-1$.  Any chamber in
the star of $v_{i}$ has the form
\begin{equation}\label{veeeye}
\sigma (v_{i},v_{k_{1}},\ldots,v_{k_{n-1}})\subset A,
\end{equation}
where $\{i, k_{1},\ldots,k_{n-1} \}\subset
\sets{n}^{\pm }$ is isotropic.  This is matched on the right side of
\eqref{homology} by the chamber
\begin{equation}\label{exp}
\sigma (v_{i},v_{k_{1}},\dots ,v_{k_{j}},x_{k_{j},k_{j+1}},\dots x_{k_{j},k_{n-1}})\subset A_{k_{j}}.
\end{equation}
In \eqref{exp}, $k_{j}$ is the first subscript reading from the left
that doesn't appear in $D_{x}$.  (The possibility that $i=k_{j}$ is
included.)  This chamber appears with the same orientation as that of
\eqref{veeeye}, and chambers of this form account for all chambers on
the left side of \eqref{homology}.

In the star of $x$, the chamber 
\begin{equation}\label{xone}
\sigma (x,x_{ij},x_{i,k_{1}},\ldots,x_{i,k_{n-2}}) \subset A_{i}
\end{equation}
is canceled by the chamber
\begin{equation}\label{xtwo}
\sigma (x,x_{ij},x_{j,k_{1}},\ldots,x_{j,k_{n-2}}) \subset A_{j},
\end{equation}
exactly as in the generic case.  If $\# D_{x}<n-1$, then every chamber
in the star of $x$ is of this form up to symmetry.  If $\# D_{x}=n-1$,
the remaining chambers in the star of $x$ have the form
\begin{equation}\label{xthree}
\sigma (x,v_{i_{1}},\dots v_{i_{n-1}}),
\end{equation}
where $D_{x} = \{i_{1},\dots i_{n-1} \}$.  This chamber will appear in
$A_{k}$ and $A_{\bar k}$ with opposite orientations, where $k$ and $\bar k$
are the unique elements that extend $D_{x}$ to an isotropic subset.

The remaining chambers on the right side of \eqref{homology} appear in
the stars of the $x_{ij}$ where $\{i,j \}\not \subset D_{x}$.  These
cancel exactly as in \eqref{one12:a}--\eqref{one12:d}.  

Finally we consider the case $\# D_{x}=n$.  This case is slightly
different, since $x$ is in the span of the $\{v_{i}\mid i\in D_{x} \}$.
Again we begin with the star of the $v_{i}$, and consider the chamber 
\[
\sigma (v_{i},v_{k_{1}},\dots v_{k_{n-1}}).
\]
If $D_{x}\not =\{i,k_{1},\dots k_{n-2} \}$, then this chamber is
matched by the chamber
\begin{equation}
\sigma (v_{i},v_{k_{1}},\dots ,v_{k_{j}},x_{k_{j},k_{j+1}},\dots, x_{k_{j},k_{n-1}})\subset A_{k_{j}}.
\end{equation}
exactly as in \eqref{exp}.  Otherwise, if $D_{x}=\{i,k_{1},\dots k_{n-2} \}$,
then this chamber is matched by 
\[
\sigma (v_{i},v_{k_{1}},\dots v_{k_{n-2}},x)\subset A_{\overline{k_{n-1}}}.
\]

Now consider the star of $x$.  If $\{i,j \}\cap D_{x}=\emptyset $,
then we have cancellation as in \eqref{xone} and \eqref{xtwo}.
Otherwise, write $D_{x}=\{i_{1},\dots i_{n} \}$.  Then in the
remaining chambers, 
\[
\sigma (x,v_{i_{1}},\dots ,v_{i_{n-1}})\subset A_{\overline{i_{n}}}
\]
cancels
\[
\sigma (x,v_{i_{1}},\dots ,v_{i_{n-2}},v_{i_{n}})\subset
A_{\overline{i_{n-1}}}.
\]

Finally, the chambers in the star of the $x_{ij}$ with $\{i,j \}\not
\subset D_{x}$ cancel exactly as in \eqref{one12:a}--\eqref{one12:d}.

\end{proof}  

\begin{remark}\label{more.general}
Theorem \ref{th1} can be proven in more generality than
stated here.  For example, the proof applies to buildings associated
to the odd orthogonal groups $SO_{2n+1}$, and can be modified to work
for buildings associated to the even orthogonal groups
$SO_{2n}$.  Using the notion of a ``perspectivity'' \cite{tits.II}, we
may prove a similar result for buildings of type $G_{2}$.  We would
like to have a ``building-theoretic'' proof of Theorem \ref{th1}.
\end{remark}

\section{Finiteness}\label{finiteness}
Throughout this section we assume that $\OK$ is a euclidean ring with
respect to the norm $\|\phantom{m}\|\colon \OK\rightarrow \Z
^{\geq 0}$.  We also assume, using Proposition \ref{sms.props} (2),
that all modular symbols have integral columns.

\subsection{}\label{introfour}
If $m\in Sp_{2n}(\OK )$, we call the class $[m]\in H_{n-1}(\partial
\bar X;\Z )$ a {\em unimodular symbol}.  The goal of this section is
to prove that the unimodular symbols span $H_{n-1}(\partial \bar X;\Z
)$. We begin with a notion that lets us measure how far a symplectic
modular symbol is from being unimodular.  In contrast to the special
linear case (\S\ref{speciallinearcase}), we use the symplectic pairing
rather than the determinant as a measure of non-unimodularity.

\begin{definition}\label{depth}
Let $[m]$ be a symplectic modular symbol with primitive columns.  The
{\em depth} of $[m] = [v_{1},\ldots,v_{n};v_{\bar n},\ldots,v_{\bar
1}]$ is the number
$$d(m):=\Max_{i\in \sets{n}}\Bigl\{\bigl\|\langle v_{i},v_{\bar
\imath}\rangle\bigr \| \Bigr\}. $$
\end{definition} 

Notice that $d(m)=1$ implies that each $\langle v_{i},v_{\bar
\imath}\rangle\in \OK^{\times }$.  Hence if $d(m)=1$ we may divide the
columns of $m$ by appropriate units to obtain $m'\in Sp_{2n}(\OK)$
satisfying $[m]=[m']$.  Thus to show that the unimodular symbols span
$H_{n-1}(\partial \bar X;\Z )$, it is sufficient to show that through a
homology we may replace a modular symbol with depth $>1$ with a cycle
of modular symbols, each of which has smaller depth.

\begin{lemma}\label{sympcandexist}
Let $[m]=[v_{1},\ldots,v_{n};v_{\bar n},\ldots,v_{\bar 1}]$ be a
symplectic modular symbol with primitive columns.  If $d(m)>1$, then
$i (v_{1},\ldots,v_{n};v_{\bar n},\ldots,v_{\bar 1})>1$, and there 
exists a candidate for the $v_{i}$.
\end{lemma}

\begin{proof}
Assume that $i(v_{1},\ldots,v_{n},v_{\bar
n},\ldots,v_{\bar 1})=1$.  Then the lattice generated by the $v_{1}$
is $\OK ^{n}$, and thus 
$m\in GL_{2n}(\OK)$.  Therefore
$$\det (m) = \prod _{i\in \sets{n}} \langle v_{i},v_{\bar \imath}\rangle \in
\OK^{\times }, $$
and $d(m)=1$, a contradiction.  This implies $i(v_{1},\ldots,v_{n},v_{\bar
n},\ldots,v_{\bar 1})>1$, and a candidate exits by Proposition
\ref{cand.exist}. 
\end{proof}

\subsection{}
Suppose that $[m]=[v_{1},\ldots,v_{n},v_{\bar
n},\ldots,v_{\bar 1}]$ is a symplectic modular symbol with $d (m)>1$,
and let $x$ be a candidate for $m$.  As in \S\ref{speciallinearcase}, write
$$x=\sum _{i\in \sets{n}^{\pm }}q_{i}v_{i} $$
with $q_{i}\in K$ satisfying $0\leq \|q_{i}\|<1$.  Define 
$[m_{i}]$ as in Definition \ref{msubi}, so that 
\begin{equation}\label{hom}
[m]=\sum _{i\in \sets{n}^{\pm }\setminus D_{x}} [m_{i}].
\end{equation}
Notice that for $i\in \sets{n}^{\pm }\setminus D_{x}$, we have 
$$\|\langle x, v_{i}\rangle \| =  \|q_{\bar \imath}\langle v_{\bar \imath}, v_{i}\rangle \|<\|\langle v_{\bar \imath}, v_{i}\rangle \|\leq d(m),$$ 
so that the norm of at least one of the symplectic pairings in each
$m_{i}$ has decreased.  However, it is not true that $d(m_{i})<d(m)$
in general, and so more than just the construction of $x$ is required
to show that the unimodular symbols span.  The following example
illustrates our strategy.

\begin{example}\label{sp6.claim}
Consider the case of $G=Sp_{6}(K)$.  As in Example \ref{sp6.example},
$[m]$ corresponds to an octahedron in $\Proj ^{5}(K)$.  Suppose that
$x$ is a candidate for $[m]$, and that $D_{x}=\emptyset $.  Construct
the modular symbols $[m_{i}]$.  Then the configuration in $\Proj
^{5}(K)$ corresponding to $[m_{1}]$ is the octahedron with vertices
$v_{1}$, $x$, and the four constructed points lying on the lines in
the original configuration in the link of $v_{1}$.  (See Figure
\ref{link.fig}.)
\begin{figure}[ht]
\begin{center}
\psfrag{v1}{$v_{1}$}
\psfrag{v3}{$v_{3}$}
\psfrag{v2}{$v_{2}$}
\psfrag{vb3}{$v_{\bar 3}$}
\psfrag{vb2}{$v_{\bar 2}$}
\psfrag{x12}{$x_{12}$}
\psfrag{x13}{$x_{13}$}
\psfrag{x1b3}{$x_{1\bar 3}$}
\psfrag{x1b2}{$x_{1\bar 2}$}
\psfrag{x}{$x$}
\includegraphics[scale = .7]{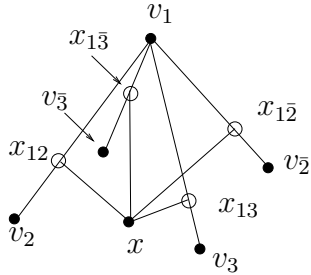}
\end{center}
\caption{\label{link.fig}The configuration corresponding to $[m_{1}]$.}
\end{figure}

The point $x$ has been chosen so that $\|\langle x, v_{1}\rangle \| <
d(m)$, but in general $\|\langle x_{12},x_{1\bar 2}\rangle\|$ and
$\|\langle x_{13}, x_{1\bar 3}\rangle\|$ will be larger than $d(m)$.
However, we claim that $i(x_{12},x_{13},x_{1\bar 3},x_{1\bar 2})>1$,
and that we may identify the quadruple $(x_{12},x_{13},x_{1\bar
3},x_{1\bar 2})$ with a modular symbol $[m']$ for $Sp_{4} (K)$.  This
means we may argue inductively and use the reduction of this ``link''
modular symbol to reduce $[m_{1}]$. (See Figure \ref{link2.fig}.)

\begin{figure}[ht]
\begin{center}
\includegraphics[scale = .7]{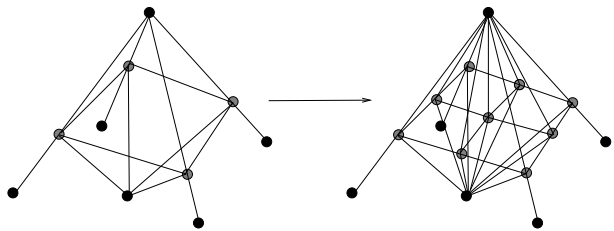}
\end{center}
\caption{\label{link2.fig}Reducing the link modular symbol.}
\end{figure}
\end{example}

\subsection{}\label{hermite}
To study the $[m_{i}]$, 
we first prove a version of Hermite normal form for our matrix representatives.

\begin{lemma}\label{hnf}
Let $m\in M_{2n}(\OK)$ be a matrix with nonzero, primitive columns,
and assume that $m$ satisfies the isotropy condition.  Then there is a
matrix $\gamma \in Sp_{2n}(\OK)$ such that $\gamma m$ is upper
triangular and $\gamma m$ satisfies the isotropy condition.
\end{lemma}

\begin{proof}
The left action of $Sp_{2n}(\OK)$ corresponds to row operations on
$m$, so we begin by listing the elementary symplectic row operations.
Let $\{r_{i}\,|\,i\in \sets{n}^{\pm } \}$ denote the rows of $m$, and
use the notation $a\gets b$ to mean that the row vector $a$ is
to be replaced by the expression $b$.  Then
we find that we may effect the following:

\begin{itemize}
\item [(T1)] $r_{i}\gets r_{i}+r_{\bar \imath}$
\item [(T2)] for $1\leq i<k\leq n$, 
\begin{align*}
r_{i}&\gets r_{i}+r_{k}\\
r_{\bar k}&\gets r_{\bar k}-r_{\bar \imath}
\end{align*}
\item [(T3)] for $1\leq i<k\leq n$, 
\begin{align*}
r_{i}&\gets r_{i}+r_{\bar k}\\
r_{k}&\gets r_{k}+r_{\bar \imath}
\end{align*}
\item [(P1)]  
\begin{align*}
r_{i}&\gets r_{\bar \imath}\\
r_{\bar \imath}&\gets -r_{i}
\end{align*}
\item [(P2)] if $\tau \in S_{n}$,
then
$$r_{i}\gets r_{\tau (i)},$$
\end{itemize}

Here {\em T} stands for {\em transvection} and {\em P} stands for {\em
permutation}.  In (P2) we mean that $\tau $ permutes the first $n$
rows among themselves, with the action on the last $n$ rows determined
by $\tau (\bar \imath ) := \overline{\tau (i)}$.  Furthermore, the
inverses and transpositions of these operations are also elementary
row operations.

Now, using operations of type (T1) and (P1), we may carry $m$ into the
following form:
$$\left(
\begin{array}{cccc} 
m_{11}&m_{12}&\cdots &m_{1\bar 1}\\
\vdots&\vdots&&\vdots\\
m_{n,1}&m_{n,2}&\cdots &m_{n,\bar 1}\\
0&m_{\bar n,2}&\cdots &m_{\bar n,\bar 1}\\
\vdots&\vdots&&\vdots\\
0&m_{\bar {1}2}&\cdots &m_{\bar 1\bar 1}
\end{array}\right). $$
In the first column, the entry in row $i$ is the greatest common
divisor of the original entries in row $i$ and row $\bar \imath$.  Next, using
operations of type (T2) and (P2), we can take $m$ into the form
$$\left(
\begin{array}{cccc} 
m_{11}&m_{12}&\cdots &m_{1\bar 1}\\
0&m_{22}&\cdots &m_{2\bar 1}\\
\vdots&\vdots&&\vdots\\
0&m_{\bar{1}2}&\cdots &m_{\bar 1\bar 1}
\end{array}\right). $$
Since left multiplication by $\gamma \in Sp_{2n}(\OK)$ preserves the
isotropy condition, it follows that $m$ actually has the form
$$\left(
\begin{array}{cccc} 
m_{11}&m_{12}&\cdots &m_{1\bar 1}\\
0&m_{22}&\cdots &m_{2\bar 1}\\
\vdots&\vdots&&\vdots\\
\vdots&m_{\bar{2}2}&\cdots &m_{\bar 2\bar 1}\\
0&\cdots&0&m_{\bar 1\bar 1}
\end{array}\right). $$
Now we may apply the induction hypothesis to the middle $(2n-2)\times
(2n-2)$ block to complete the proof of the lemma.
\end{proof}

\subsection{}
Let $[m]$ be a symplectic modular symbol with $m$ upper triangular,
and let $x$ be a candidate for $m$.  Without loss of generality, we
may assume 
$v_{1}=e_{1}$ by Proposition \ref{sms.props} (3), and that 
$1\not \in D_{x}$.  We want to study
the modular symbol $[m_{1}]$ from \eqref{hom} in greater detail.

\begin{lemma}\label{msubilemma}
Suppose that $m$ is upper triangular with $v_{1}=e_{1}$.  Let the
points 
$$\bigl\{x_{1,j}\,\bigm|\,j\in \sets{n}^{\pm }\setminus\{1,\bar 1 \}
\bigr\} $$
and the matrix $m_{1}$ be defined as in Definition \ref{msubi}.
Let $X$ be the $2n\times (2n-2)$ matrix with columns the vectors  
$$x_{12},\ldots,x_{1,n},x_{1,\bar n},\ldots,x_{1\bar 2}.$$
That is, $X$ is the matrix obtained by deleting the first and last columns
from $m_{1}$.
Define points $w_{j}$ for $j\in \sets{n}^{\pm }\setminus\{1,\bar 1 \}$ by
$$w_{j}=\langle e_{j},x\rangle e_{1} - \langle e_{1},x \rangle e_{j}, $$
and let $W$ be the $2n\times (2n-2)$ matrix with columns the vectors 
$$w_{2},\ldots,w_{n},w_{\bar n},\ldots,w_{\bar 2}. $$
Let $m'$ be the central $(2n-2)\times (2n-2)$ block of $m_{1}$.  Then
$$X = W m'.$$
\end{lemma}

\begin{proof}
Given a column vector $v$, let $v^{k}$ be the entry in the $k$-th
row.  
Suppose first that $j\leq n$.  Then computing using Definition \ref{msubi},
we find
$$x^{1}_{1,j}=x^{\bar 2}m_{2,j}+\cdots +x^{\bar \jmath}m_{jj},$$
and 
$$x^{k}_{1,j}=-x^{\bar 1}m_{kj}\quad \hbox{for $2\leq k\leq j$.} $$
On the other hand, $W$ has the form
$$\left(
\begin{array}{cccc} 
x^{\bar 2}&x^{\bar 3}&\cdots& x^{2}\\
-x^{\bar 1}&0&\cdots &0\\
0&-x^{\bar 1}&\ddots&\vdots\\
\vdots&\ddots&\ddots&0\\
0&\cdots  &0&-x^{\bar 1}
\end{array}\right),$$  
and so $x_{1,j}$ is this matrix times the $j$-th
column of $m'$.  The computation is similar if $j>n$, so the result
follows.
\end{proof}

\subsection{}
We come now to the main result of this article.  

\begin{theorem}\label{th2}
As $m$ ranges over $Sp_{2n}(\OK )$, the classes $[m]$ span
$H_{n-1}(\partial \bar X ;\Z )$.  
\end{theorem}

\begin{proof}
By the paragraph following Definition \ref{depth}, it is sufficient to
show that if $[m]$ satisfies $d (m)>1$, then $[m]=\sum [m_{\alpha }]$,
where $d (m_{\alpha })<d (m)$.  We proceed by induction.  Since
$SL_{2}(\OK )=Sp_{2}(\OK )$, the statement is true by the usual
modular symbol algorithm.  So we assume the statement is true for
$Sp_{2k}(\OK )$ with $k<n$.

Assume that $d(m)>1$, and let $x$ be a candidate for $m$.  Write 
\begin{equation}\label{finalhom}
[m]=\sum _{i\in \sets{n}^{\pm }\setminus D_{x}} [m_{i}].
\end{equation}
We permute the columns of each $[m_{i}]$ so that $v_{i}$ is the first
column.  

Choose an $[m_{i}]$ from the right of \eqref{finalhom}.
Permuting labels if necessary, we can assume $i=1$.  Since candidate
selection and the relation \eqref{homology} are $Sp_{2n}(\OK
)$-equivariant, using Lemma \ref{hnf} we may assume $m_{1}$ is
upper triangular.  As in Lemma \ref{msubilemma}, we construct the
matrix $m'$ and the vectors $w_{j}$.  Since $m'$ is a $(2n-2)\times
(2n-2)$ matrix satisfying the isotropy condition, we have that $[m']$
corresponds to a class in $H_{n-2}(\B ';\Z )$, where $\B '$ is the
building associated to $Sp_{2n-2}(K)$.  By the induction hypothesis we
may write
\begin{equation}\label{eqA}
[m']=\sum_{\alpha \in A} [m'_{\alpha }], 
\end{equation}
where $m'_{\alpha }\in Sp_{2n-2}(\OK)$, and the sum is finite.

We may use the $[m'_{\alpha }]$ to write 
\begin{equation}\label{eqB}
[m_{1}]=\sum_{\alpha \in A} [m_{\alpha }] 
\end{equation}
as follows.  Each $m'_{\alpha }$ corresponds to an endomorphism of the
$\OK $-module generated by the $w_{j}$.  We may apply $m'_{\alpha }$
to $W$ to produce a $2n\times (2n-2)$ matrix $W_{\alpha }$.  Then
$m_{\alpha }$ is the matrix with first column $e_{1}$, last column
$x$, and with middle columns $W_{\alpha }$.
The induction hypothesis asserts that the columns of $W_{\alpha }$
form an $\OK $-basis of the $\OK $-module generated by the $w_{j}$.
In particular, 
\[
d(m_{\alpha })\leq \|\langle w_{j},w_{\bar \jmath} \rangle \| =
\|\langle x, v_{1} \rangle \|^{2}.
\]
We claim that we may reduce each $[m_{\alpha }]$ further
to write
$$[m_{\alpha }] = \sum _{\beta \in B} [m_{\alpha \beta }], $$
where
$$ d(m_{\alpha\beta  })\leq \|\langle x, v_{1} \rangle \|.$$

To see this, consider the $j$ and $\bar \jmath$ columns of $W$:
\begin{align*}
w_{j}&=(x^{\bar \jmath},0,\dots,0,-x^{\bar 1},0,\dots ,0)^{t} \quad\hbox{and} \\
w_{\bar \jmath}&=(-x^{j},0,\dots,0,-x^{\bar 1},0,\dots ,0)^{t},
\end{align*}
where $-x^{\bar 1}$ appears in the $j$-th and $\bar \jmath$-th rows respectively.
Assume first that these vectors are primitive.
Then some multiple of $w_{j}$ added to $w_{\bar \jmath}$ will be
divisible by $x^{\bar 1}$.  This means $w_{j}$ and $w_{\bar \jmath}$
generate a lattice $L$ satisfying 
$$i(L) \geq \| x^{\bar 1} \| =\| \langle x, v_{1} \rangle\|.  $$
Therefore we may apply the $SL_{2}$-modular symbol algorithm to the
pair $(w_{j}, w_{\bar \jmath})$ to reduce $L$ to $(L\otimes K) \cap \OK
^{n}$.  If either $w_{j}$ or $w_{\bar \jmath}$ is not primitive, then
a simple modification of this argument achieves the same end.

Applying this to each pair $(w_{j}, w_{\bar \jmath})$, and adapting 
the construction that passes from \eqref{eqA} to \eqref{eqB}, we find 
\begin{equation}\label{eqC}
[m_{\alpha }] = \sum _{\beta \in B} [m_{\alpha \beta }], 
\end{equation}
with
$$ d(m_{\alpha\beta  })\leq \|\langle x, v_{1} \rangle \|.$$
Together \eqref{eqB} and \eqref{eqC} imply 
\[
[m_{1}] = \sum_{\substack{\alpha \in A\\
\beta \in B}} [m_{\alpha\beta  }]\quad \hbox{with
$d(m_{\alpha\beta  }) \leq \langle x, v_{1} \rangle<d (m)$.} 
\]
Since this argument may be applied to any of the $[m_{i}]$ from
\eqref{finalhom}, this completes the proof of the theorem.
\end{proof}

\begin{corollary}\label{finspsetsymp}
If $\Gamma \subset Sp_{2n}(\OK)$ is torsion-free of finite index, and
$N$ is the cohomological dimension of $\Gamma $, then the classes
$[m]^{*}_\Gamma $ provide a finite spanning set of $H^{N}(\Gamma ;\Z
)$.
\end{corollary}
\bibliographystyle{amsplain}
\bibliography{sympms}

\end{document}